\documentclass[10pt,a4paper]{amsart}

\pdfoutput=1

\usepackage{amsmath}
\usepackage{amssymb}
\usepackage{amsopn}
\usepackage{epsfig}
\usepackage{amsfonts}
\usepackage{latexsym}

\usepackage{subfigure}

\title[Bounds for Diophantine approximation]{Sharp bounds for symmetric and asymmetric Diophantine approximation}
\author{Cor Kraaikamp} 
\address{Delft University of Technology and Thomas Stieltjes Institute for
Mathematics, DIAM, Mekelweg 4, 2628 CD Delft, the Netherlands}
\email{c.kraaikamp@tudelft.nl}
\author{Ionica Smeets}
\address{Universiteit Leiden and Thomas Stieltjes Institute for
Mathematics, Niels Bohrweg 1, 2333 CA Leiden, the Netherlands}
\email{smeets@math.leidenuniv.nl}
\date{\today}


\usepackage[pagebackref=true]{hyperref}
\hypersetup{   pdfborder={0 0 0},   
              colorlinks = false}

\renewcommand*{\backref}[1]{}
\renewcommand*{\backrefalt}[4]{%
  \ifcase #1 %
    No citations.
  \or
    (Page #4).%
  \else
    (Pages #4).%
  \fi%
}

\newcommand{\R}{\mathbb R}
\newcommand{\Z}{\mathbb Z}
\newcommand{\Q}{\mathbb Q}
\newcommand{\N}{\mathbb N}

\newtheorem{Theorem}{Theorem}[section]
\newtheorem{Proposition}[Theorem]{Proposition}
\newtheorem{Lemma}[Theorem]{Lemma}

\newtheorem{Corollary}[Theorem]{Corollary}

\theoremstyle{definition}

\theoremstyle{remark}
\newtheorem{Rmk}[Theorem]{Remark}
\newtheorem{Ex}[Theorem]{Example}

\begin{document}
\begin{abstract}
In 2004, J.C. Tong found bounds for the approximation quality of a regular continued fraction convergent of a rational number, expressed in bounds for both the previous and next approximation.
We sharpen his results with a geometric method and give both sharp upper and lower bounds. We also calculate the asymptotic frequency that these bounds occur.
\end{abstract}

\setlength{\parindent}{0pt}
\setlength{\parskip}{1.8ex}

\maketitle

\section{Introduction}
In 1891, Hurwitz showed in~\cite{H} that for every irrational number $x$ there exist infinitely many co-prime integers $p$ and $q$, with $q >0$, such that
$$
\left| x-\frac{p}{q}\right| < \frac{1}{\sqrt{5}} \frac{1}{q^2},
$$
where the constant $1/\sqrt{5}$ is ``best possible,'' in the sense that it cannot be replaced by a smaller constant.\smallskip\

Let $x$ be a real irrational number, with regular continued fraction (RCF) expansion
\begin{equation}\label{rcf}
x=a_0+\frac{\displaystyle 1}{\displaystyle a_1+\frac{1}{a_2+\ddots
+\displaystyle \frac{1}{a_n+\ddots}}} = [a_0; a_1,a_2,\dots,
a_n,\dots ].
\end{equation}
Here $a_0\in\Z$ is such, that $x-a_0\in [0,1)$, and $a_n\in  \N$ for $n\geq 1$. Finite truncation in~(\ref{rcf}) yields the convergents $p_n/q_n$, $n\geq 0$, i.e.,
$$
\frac{p_n}{q_n}= [a_0; a_1,a_2,\dots, a_n],\quad \text{for $n\geq
1$}.
$$
The partial coefficients $a_{n}$ can be found from the regular continued fraction map\ $T\, :\, [0,1)\rightarrow
[0,1)$, defined by
$$
T(x) := \left\{ \frac{1}{x}\right\} = \frac{1}{x}-\left\lfloor \frac{1}{x}\right\rfloor ,\, x\ne 0;\quad T(0) := 0,
$$
where $\lfloor x \rfloor$ denotes the largest integer smaller than or equal to $x$.

Borel showed in~\cite{B} that for all $n\geq 1$,
\begin{equation}\label{borel}
\min \{ \Theta_{n-1},\Theta_n,\Theta_{n+1}\} <
\frac{1}{\sqrt{5}},
\end{equation}
where the \emph{approximation coefficients} $\Theta_n$ of $x$ are defined by
\begin{equation}\label{theta-n}
\Theta_n=\Theta_n(x)=q_n^2\left| x-\frac{p_n}{q_n}\right| ,\qquad
\text{for $n\geq 0$}.
\end{equation}

Hurwitz' result is a direct consequence of Borel's result, and a classical theorem by Legendre, which states that if $p$ and $q$ are two co-prime integers with $q>0$, satisfying
$$
\left| x-\frac{p}{q}\right| < \frac{1}{2} \frac{1}{q^2},
$$
then there exists an $n\in\N$, such that $p=p_n$ and $q=q_n$.\smallskip\

%

Over the last century Borel's result~(\ref{borel}) has been refined in various ways. For example, in~\cite{Fu}, \cite{Mu},
and~\cite{BMcL}, it was shown that
$$
\min \{ \Theta_{n-1},\Theta_n,\Theta_{n+1}\} <
\frac{1}{\sqrt{a_{n+1}^2+4}},\qquad \text{for $n\geq 0$},
$$
while J.C. Tong showed in~\cite{Tong} that the ``conjugate property'' holds
$$
\max \{ \Theta_{n-1},\Theta_n,\Theta_{n+1}\} >
\frac{1}{\sqrt{a_{n+1}^2+4}},\qquad \text{for $n\geq 0$}.
$$
Also various other results on Diophantine approximation  have been obtained, starting with Dirichlet's observation from~\cite{Di}, that
$$
\left| x-\frac{p_n}{q_n}\right| <  \frac{1}{q_nq_{n+1}},\qquad
\text{for $n\geq 0$},
$$
which lead to various results in symmetric and asymmetric Diophantine approximation; see e.g.~\cite{Tjnt88},
\cite{Tsnf89}, \cite{Kindag90}, and~\cite{Kjnt94}.\medskip\

Define for $x$ irrational the number $C_n$ by
\begin{equation}\label{chap1 approximationcoefficientsTong}
x-\frac{p_n}{q_n}=\frac{(-1)^n}{C_nq_nq_{n+1}},\qquad \text{for
$n\geq 0$}.
\end{equation}
Tong derived in~\cite{Tsnf89} and~\cite{Tams91} various properties of the sequence $(C_n)_{n\geq 0}$, and of the related sequence
$(D_n)_{n\geq 0}$, where
\begin{equation}
\label{def: chap1 D_n}
D_n=[a_{n+1};a_n,\dots,a_1]\cdot [a_{n+2};a_{n+3},\dots ] = \frac{1}{C_n -1} ,\qquad
\text{for $n\geq 0$}.
\end{equation}
%

Recently, Tong~\cite{Tca04} obtained  the following theorem, which covers many previous results.
\begin{Theorem}\label{thm:tong2004} {\rm(Tong)}
Let $x=[\, a_0; a_1,\, a_2,  \ldots ,\, a_n, \ldots\, ]$ be an irrational number. If $r>1$ and $R>1$ are two real numbers and
$$
\begin{aligned}
M_{\rm{Tong}} = \quad &\frac{1}{2} \left( \frac{1}{r} + \frac{1}{R} + a_n a_{n+1}\left(1 + \frac{1}{r} \right)\left(1 + \frac{1}{R}\right)\right.\nonumber\\
&+ \left.\sqrt{\left[\frac{1}{r} + \frac{1}{R} + a_n a_{n+1}\left(1 + \frac{1}{r} \right)\left(1 +  \frac{1}{R}\right)\right]^2 - \frac{4}{rR}}\,\right), \label{eq: MTong}
\end{aligned}
$$
then
\begin{enumerate}
 \item $D_{n-2}<r$ and $D_n<R$ imply $D_{n-1}>M_{\rm{Tong}}$;
\item $D_{n-2}>r$ and $D_n>R$ imply $D_{n-1}<M_{\rm{Tong}}$.
\end{enumerate}
\end{Theorem}
Tong derived a similar result for the sequence $C_{n}$, but it is incorrect. We state this result, give a counterexample and present a correct version of it in Section~\ref{sec: C}.  

In Section~\ref{sec:smaller} we prove the following result.
\begin{Theorem}
\label{th: intro smaller}
Let $r,R>1$ be reals and put
$$
F =  \frac{r(a_{n+1}+1)}{a_{n}(a_{n+1}+1)(r+1)+1}  \quad \textrm{and} \quad G =   \frac{R(a_{n}+1)}{(a_{n}+1)a_{n+1}(R+1)+1} .
$$

If $D_{n-2}<r$ and $D_{n}<R$, then 
\begin{enumerate}
  \item If $\displaystyle{r-a_{n} \geq G}$ and $\displaystyle{R-a_{n+1} <F}$, then  
\[
D_{n-1}>\frac{a_{n+1}+1}{R-a_{n+1}},
\]
 \item If $\displaystyle{r-a_{n} < G \textrm{ and } R-a_{n+1} \geq F}$, 
then  
\[
D_{n-1}> \frac{a_{n}+1}{r-a_n} ,
\]
 \item In all other cases
\[
D_{n-1}>M_{\rm{Tong}}.
\]
\end{enumerate}
These bounds are sharp.
\end{Theorem}

The outline of this paper is as follows. We derive elementary properties of the sequence $D_n$ in Section~\ref{sec: tong natural ext}. In Section~\ref{sec:smaller} we prove Theorem~\ref{th: intro smaller} that gives a sharp lower bound for the minimum of $D_{n-1}$ in case $D_{n-2}<r$ and $D_{n}<R$ for real numbers  $r,R>1$. We prove a similar theorem for the case that $D_{n-2}>r$ and $D_n>R$ in Section~\ref{sec:larger}. In Section~\ref{sec: probabilities} we calculate the asymptotic frequency that simultaneously $D_{n-2}>r$ and $D_{n}>R$. Finally we correct Tong's result for $C_{n}$ in Section~\ref{sec: C} and give the sharp bound in this case. 


%

\section{The natural extension}
\label{sec: tong natural ext}
Define the space $\Omega = [0,1) \times [0,1]$ and define $\mathcal{T}:\Omega \rightarrow \Omega$ as
\[\mathcal{T}(x,y)=\left(\frac{1}{x}-\left\lfloor \frac{1}{x}\right\rfloor,\frac{1}{a_1(x)+y}\right).\]

The following theorem was obtained in 1977 by Nakada \emph{et al.} \cite{NIT}; see also~\cite{Nakada} and \cite{IK}.

\begin{Theorem}
\label{th: ergodic}
Let $\nu$ be the probability measure on $\Omega$ with density $d(x,y)$, given by
\begin{equation}
\label{eq: density}
d(x,y)=\frac{1}{\log 2}\frac{1}{(1+xy)^2}, \quad (x,y) \in \Omega;
\end{equation}
then $\nu$ is the invariant measure for $\mathcal{T}$. Furthermore, the dynamical system $(\Omega,\nu,\mathcal{T})$ is an ergodic system.
\end{Theorem}

The system $(\Omega,\nu,\mathcal{T})$ is the natural extension of the ergodic dynamical system $([0,1),\mu,T)$, where $\mu$ is the so-called Gauss-measure, the probability measure on $[0,1)$ with density
\[d(x)=\frac{1}{\log 2} \frac{1}{1+x}, \quad x \in [0,1).\]

This natural extension plays a key role in the proofs of various results in this paper.

We write $t_n$ and $v_n$ for the ``future'' and ``past'' of $\frac{p_n}{q_n}$, respectively,
\begin{equation}
\label{eq: into tn and vn}
t_n = [a_{n+1},a_{n+2},\dots]\quad  \textrm{ and } \quad  v_n = [a_n,\dots,a_1].
\end{equation}

Furthermore, $t_0=x$ and $v_0 =0$. 

The approximation coefficients may be written in terms of $t_{n}$ and $v_{n}$
$$
\Theta_n\, =\, \frac{t_n}{1+t_nv_n}\, \qquad \textrm{and} \qquad  \Theta_{n-1}\, =\, \frac{v_n}{1+t_nv_n}\, ,\qquad n\geq 1.
$$

\begin{Lemma}
Let $x=[a_0;a_1,a_2,\dots]$ be in $\R \setminus \Q$ and $n \geq 2$ be an integer. The variables $D_{n-2}, D_{n-1}$ and $D_n$ can be expressed in terms of future $t_n$, past $v_n$ and digits $a_n$ and $a_{n+1}$  by
\begin{eqnarray}
\label{eq: dn-2 in tn and vn}
D_{n-2} = & D_{n-2} (t_n,v_n)&=\frac{(a_n+t_n)v_n}{1-a_nv_n},\\
\label{eq: dn-1 in tn and vn}
 D_{n-1}=&D_{n-1}(t_n,v_n)&= \frac{1}{t_n v_n}, \quad \textrm{ and}\\
 \label{eq: dn in tn and vn}
 D_n =&D_n(t_n,v_n)&=  \frac{(a_{n+1}+v_n)t_n}{1-a_{n+1}t_n}.
\end{eqnarray}
\end{Lemma}
\begin{proof}
The expression for $D_{n-1}$ follows from the definition in~(\ref{def: chap1 D_n}).
$$
\begin{aligned}
D_{n-1} &= [\,a_n;a_{n-1},\dots,a_1\,] [\,a_{n+1};a_{n+2},\dots \,]\\
&= \frac{1}{[\,0;a_n,a_{n-1},\dots,a_1\,][\,0;a_{n+1},a_{n+2},\dots \,]}=  \frac{1}{v_n t_n}.
\end{aligned}
$$
It follows in a similar way that $D_n= \frac{1}{t_{n+1}}\frac{1}{v_{n+1}}$ and using
$$
\begin{aligned}
t_{n+1} &= \frac{1}{t_n}-a_{n+1}\\
v_{n+1}&= \frac{q_n}{q_{n+1}} = \frac{q_n}{a_{n+1}q_n+q_{n-1}} = \frac{1}{a_{n+1}+v_n}
\end{aligned}
$$
we find~\eqref{eq: dn in tn and vn}. The formula for $D_{n-2}$ can be derived in a similar way.
\end{proof}

\begin{Rmk}
Of course, $D_{n-2},D_{n-1}$ and $D_n$ also depend on $x$, but we suppress this dependence in our notation.
\end{Rmk}

The following result on the distribution of the sequence $(t_{n},v_{n})_{n \geq 0}$ is a consequence of the Ergodic Theorem, and was originally obtained by W. Bosma \emph{et al.} in~\cite{BJW}, see also Chapter 4 of~\cite{IK}.

\begin{Theorem}
For almost all $x \in [0,1)$ the two-dimensional sequence
\[(t_{n},v_{n})=\mathcal{T}^{n}(x,0), \quad n\geq 0,\]
is distributed over $\Omega$ according to the density-function $d(t,v)$, as given in~(\ref{eq: density}).
\end{Theorem}
 
 Consequently, for any Borel measurable set $B \subset \Omega$ with a boundary of Lebesque measure zero, one has that
 \begin{equation}
 \label{eq: limit prob}
\lim_{n \rightarrow \infty} \frac{1}{n}\sum_{k=0}^{n-1} I_{B}(t_{n},v_{n}) = \nu(B),  
\end{equation}
 where $I_{B}$ is the indicator function of $B$. We use this result to derive the following proposition.
 
 \begin{Proposition}
 For almost all $x \in [0,1)$, and for all $R\geq 1$, the limit
 \[\lim_{n \rightarrow \infty} \frac{1}{n} \# \{1 \leq j \leq n \, | \,D_{j}(x) \leq R\}\]
exists, and equals
\begin{equation}
\label{eq: H(D)}
H(R)=1-\frac{1}{\log 2}\left(\log \left(\frac{R+1}{R}\right)+\frac{\log R}{R+1} \right).
\end{equation} 
 Consequently, for almost all $x \in [0,1)$ one has that
 \[
 \lim_{n \rightarrow \infty} \frac{1}{n} \sum_{k=0}^{n-1}D_{n}(x)=\infty.
 \]
  \end{Proposition}
 \begin{proof}
By~(\ref{eq: dn-1 in tn and vn}) and~(\ref{eq: limit prob}) for almost every $x$ the asymptotic frequency that $D_{n-1}\leq R$ is given by the measure of those points $(t,v)$ in $\Omega$ with $\frac{1}{t v} \leq R$. 
This measure equals
\[
\frac{1}{\log 2}\int_{t=\frac{1}{R}}^{1} \int_{v=\frac{1}{Rt}}^{1} \frac{{\rm d}v \, {\rm d}t }{(1+tv)^{2}} ;
\]
also see Figure~\ref{fig: expectation}. 

\begin{figure}[!ht]
\makebox[4cm][c]{\includegraphics[width=5cm]{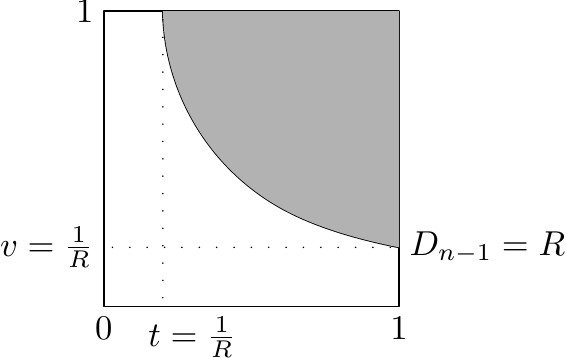}}
\caption{The curve $\frac{1}{tv}=R$ on $\Omega$. For $(t_n,v_n)$ in the gray part it holds that $D_{n-1}\leq R.$}
\label{fig: expectation}
\end{figure}

It follows that
\[
H(R) = \frac{1}{\log 2} \int_\frac{1}{R}^1 \left[ \frac{v}{1+tv} \right]_{\frac{1}{Rt}}^1 {\rm d}t = \frac{1}{\log 2} \left[ \log 2    - \log \frac{R+1}{R}  - \frac{1}{R+1}\log R \right],
\]
which may be rewritten as~(\ref{eq: H(D)}).

To calculate the expectation of $D_{n}$ we use that the density function of $D_n$ is given by $h(x)=H'(x)$, so 
\[
h(x) = \frac{1}{\log 2} \frac{\log x}{(x+1)^{2}}, \quad \textrm{for } x \geq 1.
\]
We can now easily calculate the expected value of $D_n$
\[
\lim_{n \rightarrow \infty} \frac{1}{n} \sum_{j=0}^{n-1}D_{j}(x) = \int_{1}^{\infty} x\,h(x) \, {\rm d}x = \lim_{t \rightarrow \infty }\int_{1}^{t} \frac{1}{\log 2} \frac{x \log x}{(x+1)^{2}} \, {\rm d}x =\infty.
\]
\end{proof}

 Apart for proving metric results on the $D_{n}$'s, the natural extension $(\Omega,\nu,\mathcal{T})$ is also very handy to obtain various Borel-type results on the $D_{n}$'s. 

For $a,b \in \N$ consider the rectangle $\Delta_{a,b}=\displaystyle{\left[\frac{1}{b+1},\frac{1}{b}\right) \times \left[\frac{1}{a+1},\frac{1}{a}\right)} \subset \Omega$. On this rectangle we have $a_{n}=a$ and $a_{n+1}=b$. So $(t_{n},v_{n})\in \Delta_{a,b}$ if and only if $a_{n}=a$ and $a_{n+1}=b$ . We use $a$ and $b$ as abbreviation for $a_{n}$ and $a_{n+1}$, respectively, if we are working in such a rectangle. 

We define two functions from $\left[ \frac{1}{b+1},\frac{1}{b}\right)$ to $\R$,
\begin{equation}
\label{eq: f and g}
 f_{a,r}(t)=\frac{r}{a(r+1)+t} \quad \textrm{ and } \quad  g_{b,R}(t)=\frac{R}{t}-b(R+1).
\end{equation}
From~(\ref{eq: dn-2 in tn and vn})~and~(\ref{eq: dn in tn and vn}) it follows for $(t_{n},v_{n}) \in \Delta_{a,b}$ that 
\begin{eqnarray*}
D_{n-2}<r  \quad &\textrm{ if and only if }& \quad v_n < f_{a,r}(t_n),\\ 
D_n < R \quad &\textrm{ if and only if }& \quad v_n < g_{b,R}(t_n).
\end{eqnarray*}

\begin{figure}[!ht]
\makebox[10cm][c]{\includegraphics[width=10cm]{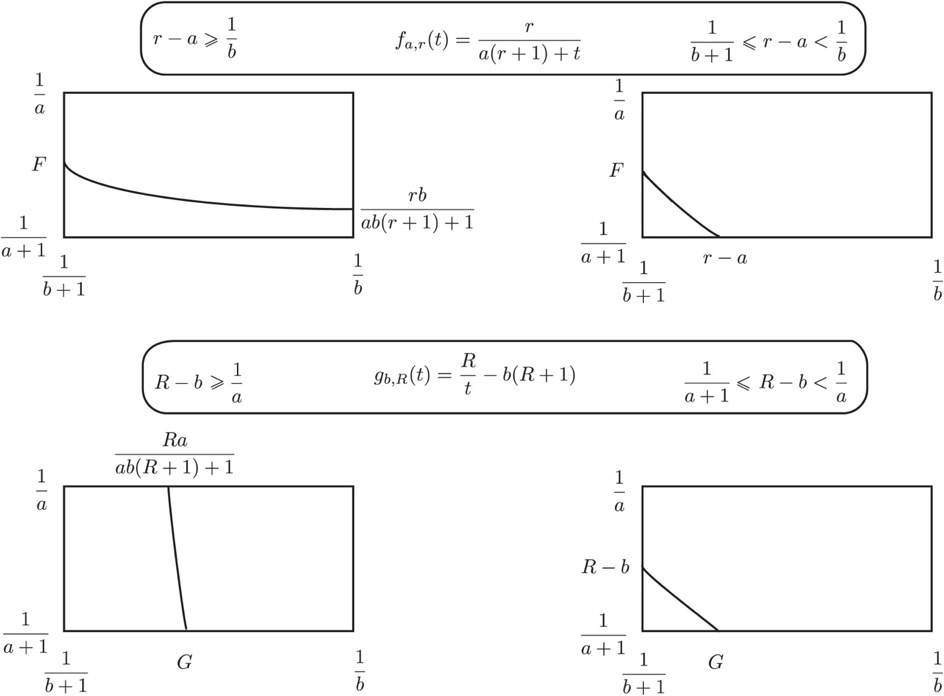}}
\caption{The possible intersection points of the graphs of $f_{a,r}$ and $g_{b,R}$ and the boundary of the rectangle  $\Delta_{a,b}$, where $a_n=a$ and $a_{n+1}=b$.}
\label{fig: cutting points}
\end{figure}

We introduce the following notation 
\begin{equation}
\label{def: F and G}
F =  \frac{r(b+1)}{a(b+1)(r+1)+1} \quad \textrm { and } \quad G =  \frac{R(a+1)}{(a+1)b(R+1)+1}.
\end{equation}
We have that $F=f_{a,r}\left(\frac{1}{b+1}\right)$ and $g_{b,R}(G)=\frac{1}{a+1}$; also see Figure~\ref{fig: cutting points}.

\begin{Rmk}
\label{rmk: strips}
The position of the graph of $f_{a,r}$ in $\Delta_{a,b}$ depends on $a$ and $r$. Obviously we always have $ f_{a,r}\left(\frac{1}{b}\right) <  f_{a,r}\left(\frac{1}{b+1}\right) = F < \frac{1}{a} $. Furthermore
$$
\begin{aligned}
f_{a,r} \left( \frac{1}{b+1}\right) \geq \frac{1}{a+1} &\quad \textrm{if and only if } \quad  r \geq a+\frac{1}{b+1},\\
f_{a,r} \left( \frac{1}{b}\right) \geq \frac{1}{a+1} &\quad \textrm{if and only if } \quad r \geq a+\frac{1}{b}.\\
\end{aligned}
$$
Similarly, the position of the graph of $g_{b,R}$ in $\Delta_{a,b}$ depends on $b$ and $R$. We always have $G < \frac{1}{b} $. Furthermore
$$
\begin{aligned}
G \geq \frac{1}{b+1} &\quad \textrm{if and only if } \quad R \geq b+ \frac{1}{a+1},\\
g_{b,R} \left( \frac{1}{b+1}\right) < \frac{1}{a} &\quad \textrm{if and only if } \quad  R < b+\frac{1}{a},\\
g_{b,R} \left( \frac{1}{b+1}\right) \geq \frac{1}{a+1} &\quad \textrm{if and only if } \quad R \geq b+\frac{1}{a+1}.\\
\end{aligned}
$$
Compare with~Figure~\ref{fig: cutting points}.

\end{Rmk}

We use the following lemma to determine where $D_{n-1}$ attains it extreme values.
\begin{Lemma}
\label{lem: increasing decreasing}
Let $a,b \in \N$, and let $D_{n-1}(t,v)=\frac{1}{tv}$ for points $(t,v) \in (0,1] \times (0,1]$. 
\begin{enumerate}
\item When $t$ is constant, $D_{n-1}$ is monotonically decreasing as a function of $v$.
\item When $v$ is constant, $D_{n-1}$ is monotonically decreasing as a function of $t$.
\item $D_{n-1}(t,v)$ is monotonically decreasing as a function of $t$ on the graph of $f_{a,r}$.
\item  $D_{n-1}(t,v)$ is monotonically increasing as a function of $t$ on the graph of $g_{b,R}$.
 \end{enumerate}
\end{Lemma}
\begin{proof}
The first two statements follow from the trivial observation
\begin{equation}
\label{eq: partial D n-1 }
\frac{\partial D_{n-1}}{\partial t}<0 \quad \textrm{ and } \quad \frac{\partial D_{n-1}}{\partial v}<0.
\end{equation}

For points $(t,v)$ on the graph of $f_{a,r}$ we find $D_{n-1}(t,v) = \frac{a(r+1)+t}{rt}$ and 
$$
\frac{\partial D_{n-1}}{\partial t} = \frac{-a(r+1)}{rt^2}<0, 
$$
which proves (3).

Finally, for points $(t,v)$ on the graph of $g_{b,R}$ we find \mbox{$D_{n-1}(t,v) = \frac{1}{R-b(R+1)t}$}. So $\frac{\partial D_{n-1}}{\partial t} >0$ and (4) is proven.
\end{proof}

\begin{Corollary}
On $\Delta_{a,b}$ the infimum  of $D_{n-1}$ is attained in the upper right corner  and its maximum in the lower left corner. To be more precise
\[
ab < D_{n-1} \leq (a+1)(b+1).
\]
\end{Corollary}

\begin{Lemma}
\label{lem: intersection f and g}
Let $a,b \in \N$, $r,R>1$, and set
\[
L=ab(r+1)(R+1), \, w=\sqrt{4LR+(r-R+L)^2)} \textrm { and } s= \left(\frac{-L+R-r+w}{2b(R+1)}\right).
\] 
On $\R_+$ the graphs of $f_{a,r}$ and $g_{b,R}$ have one intersection point, which is given by
$$
(S,f_{a,r}(S) = \left(\frac{-L+R-r+w}{2b(R+1)},\frac{2br(R+1)}{L+R-r+w}\right),
$$
The corresponding value for $D_{n-1}$ in this point is given by $M_{\rm{Tong}}$ as defined in~{\rm(\ref{eq: MTong})}.
For $x<S$ one has that $f_{a,r}(x)<g_{b,R}(x)$, while $f_{a,r}(x)>g_{b,R}(x)$ if $x>S$.
\end{Lemma}
\begin{proof}
Solving $$\frac{r}{a(r+1)+t} =\frac{R}{t}-b(R+1)$$ yields
$$
S= \frac{-L+R-r+w}{2b(R+1)}  \quad \textrm{or} \quad S= \frac{-L+R-r-w}{2b(R+1)}.
$$
Since $L > R$ the second solution is always negative, so this solution can not be in $\Delta_{a,b}$. The second coordinate follows from substituting $S= \frac{-L+R-r+w}{2b(R+1)} $ in $f_{a,r}(t)$ or $g_{b,R}(t)$.

The corresponding value for $D_{n-1}$ in this point is given by
 $$
\begin{aligned}
&D_{n-1} \left(\frac{-L+R-r+w}{2b(R+1)},\frac{2br(R+1)}{L+R-r+w}\right)=  \frac{-L-R+r-w}{r(L-R+r-w)} \\
&= \frac{-L^2+r^2-2Rr+R^2-2Lw-w^2}{r((L-R+r)^2-w^2)} = \frac{-2L^2-2Lw-2Lr-2LR}{-4RrL} \\
& =\frac12 \left(\frac{1}{r} + \frac{1}{R} + \frac{L}{Rr} + \frac{w}{Rr}\right) = M_{\rm{Tong}}.
\end{aligned}
$$
Since $\lim_{x \downarrow 0}f_{a,r}(x) = \frac{r}{a(r+1)}$ and $\lim_{x \downarrow 0}g_{b,R}(x) = \infty$, we immediately have that $f_{a,r}(x) <g_{b,R}(x)$ if $x<S$. And because there is only one intersection point on $\R_+$, it follows that  $f_{a,r}(x) > g_{b,R}(x)$ if $x>S$.
\end{proof}

\begin{figure}[!ht]

\subfigure[In case (i) and (ii) we have $r-a \geq G$ and $R-b<F$. It is allowed that $R-b < \frac{1}{a+1}$. In case (i) we have $r-a > \frac{1}{b}$ and in case (ii) $ r-a\leq \frac{1}{b} $. ]{
\hskip20mm
\includegraphics{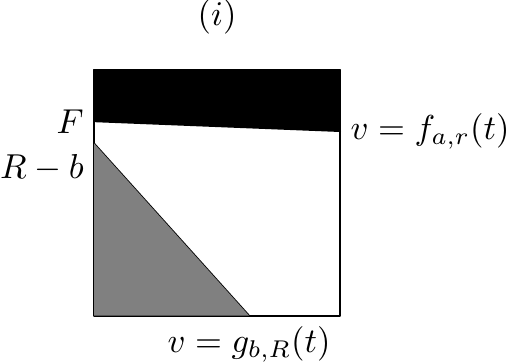}
\hskip10mm
\includegraphics{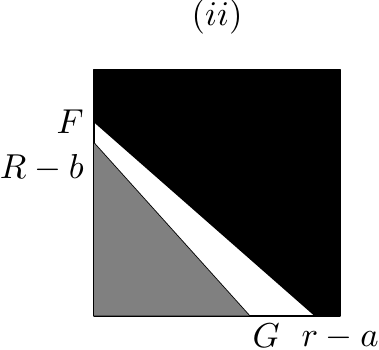}
}
\subfigure[In cases (iii) and (iv) we have $r-a < G$  and $R-b \geq F$.  It is allowed that $r-a < \frac{1}{b+1}$. In case (iii) we have $R-b > \frac{1}{a}$ and in case (iv) $ R-b \leq \frac{1}{a} $.  ]{
\includegraphics{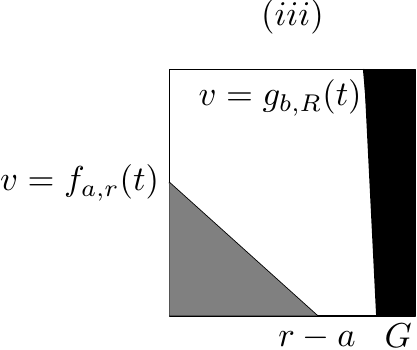}
\hskip20mm
\includegraphics{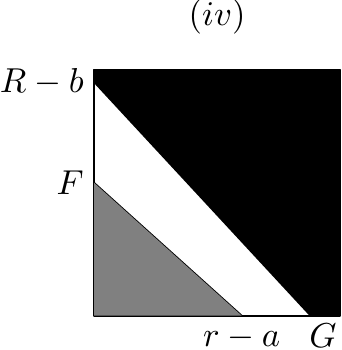}
}
\subfigure[In case (v) we have $F<\frac{1}{a+1}$ and $G<\frac{1}{b+1}$. Case (vi) contains all other cases, it can be separated in four subcases, see Figure~\ref{fig: msix}.  ]{
\hskip10mm
\includegraphics{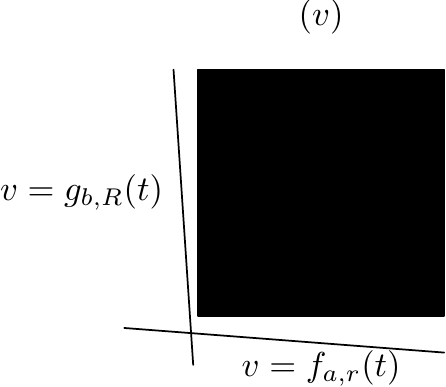}
\hskip30mm
\includegraphics{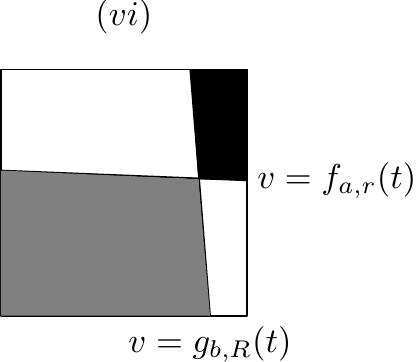}
}
\caption{The possible configurations of the graphs of $f_{a,r}$ and $g_{b,R}$ on $\Delta_{a,b}$. On the grey parts  $D_{n-2}<r$ and $D_n<R$, on the black parts  $D_{n-2}>r$ and $D_n>R$
}
\label{im: positions}
\end{figure}

\begin{Rmk}
In view of Remark~\ref{rmk: strips} and the last statement of Lemma~\ref{lem: intersection f and g} the only possible configurations for $f_{a,r}$ and $g_{b,R}$ in $\Delta_{a,b}$ are given in Figure~\ref{im: positions}.
\end{Rmk}

%


\section{The case $D_{n-2}<r$ and $D_n<R$}
\label{sec:smaller}
We assume that both $D_{n-2}$ and $D_n$ are smaller than some given reals $r$ and $R$. We prove Theorem~\ref{th: intro smaller} from the Introduction.


\begin{proof}
We consider the closure of the region containing all points $(t,v)$ in  $\Delta_{a,b}$ with $D_{n-2}(t,v)<r$ and $D_{n}(t,v)<R$. In Figure~\ref{im: positions} we show all possible configurations of this region. 

From~\eqref{eq: partial D n-1 } it follows that the extremum of $D_{n+1}$ is attained in a boundary point. Lemma~\ref{lem: increasing decreasing} implies that we only need to consider the following three points
\begin{enumerate}
\item The intersection point of the graph of $g_{b,R}$ and the line $t = \frac{1}{b+1}$, given by $\left(\frac{1}{b+1},R-b\right)$.
\item The intersection point of the graph of $f_{a,r}$ and the line $v = \frac{1}{a+1}$, given by $\left(r-a,\frac{1}{a+1}\right)$.
\item The intersection point of the graphs of $f_{a,r}$ and $g_{b,R}$, given by $M_{\rm{Tong}}$. 
\end{enumerate}

Assume $\displaystyle{r-a \geq G}$ and $R-b <F$. We know from Lemma~\ref{lem: intersection f and g} that the graphs of $f_{a,r}$ and $g_{b,R}$ can not intersect more than once in $\Delta_{a,b}$, thus we are in case (1); see Figure~\ref{im: positions} ($i$) and ($ii$). In this case the minimum of $D_{n-1}$ is given by $D_{n-1}\left(\frac{1}{b+1},R-b\right)=\frac{b+1}{R-b}$. The intersection point $(S,f_{a,r}(S))$ lies to the left of $\left(\frac{1}{b+1},R-b\right)$ and from Lemma~\ref{lem: increasing decreasing} we know that $D_{n-1}$ increases on the graph of $g_{a,r}$. We conclude that $M_{\rm{Tong}}=D_{n-1}(S,f_{a,r}(S))$ is smaller than $\frac{b+1}{R-b}$. 

Assume $\displaystyle{ r-a < G \textrm{ and } R-b \geq F}$, then we are in case (2);  see Figure~\ref{im: positions} ($iii$) and ($iv$) and the minimum is given by  $D_{n-1}=\left(r-a,\frac{1}{a+1}\right) = \frac{a+1}{r-a}$. A similar argument as before shows $M_{\rm{Tong}} < \frac{a+1}{r-a}$.

Otherwise, still assuming there are points $(t,v) \in \Delta_{a,b}$ with $D_{n-2}(t,v)<r$ and $D_{n}(t,v)<R$, we must be in case (3); see Figure~\ref{im: positions} ($vi$). The minimum follows from Lemma~\ref{lem: intersection f and g}.





These bounds are sharp since the minimum is attained in the extreme point.
\end{proof}

\begin{Ex}
Take $r=2.9$ and $R=3.6$; see Figure~\ref{im: example smaller}.

\begin{figure}[!ht]
\includegraphics[height=6cm]{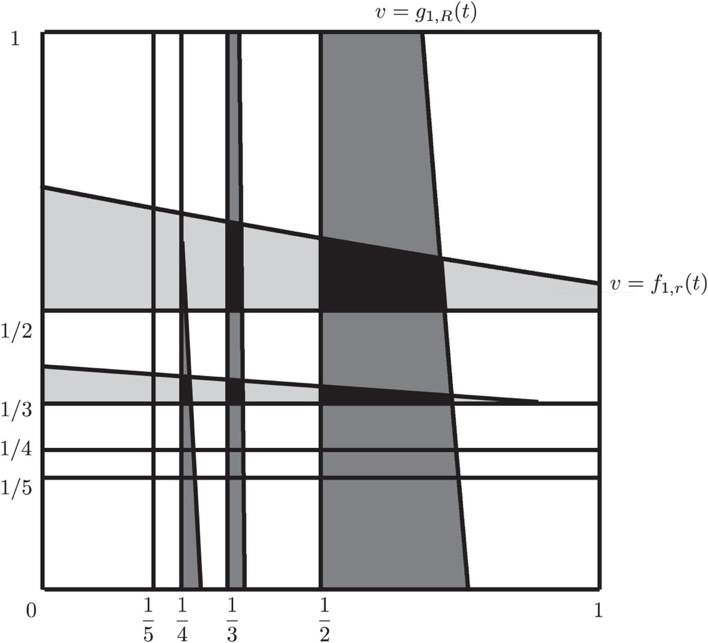}
\caption{Example with $r=2.9$ and $R=3.6$. The regions where $D_{n-2}<2.9$ are light grey, the regions where $D_n<3.6$ are dark grey. The intersection where both $D_{n-2}<2.9$ and $D_n<3.6$ is black. The horizontal and vertical black lines are drawn to identify the strips and have no meaning for the value of $D_{n-2}$ and $D_n$.}
\label{im: example smaller}\end{figure}

If $a_{n}=a_{n+1}=1$, then $r-a_{n}=1.9, R-a_{n+1}=2.6, F\approx 0.66$ and $G\approx 0.71$. Since $R-a_{n+1}>F$ we do not have case ($i$) of Theorem~\ref{th: intro smaller}. Since $r-a_{n}>G$ we are not in case ($ii$) either. So in this case $D_{n-1}>M_{\rm{Tong}}\approx 2.30$.
For the following combinations the minimum is also given by $M_{\rm{Tong}}$
\begin{eqnarray*}
a_{n}=1 \textrm{ and } a_{n+1}=2: D_{n-1} >M_{\rm{Tong}}&\approx& 4.04.\\
a_{n}=2 \textrm{ and } a_{n+1}=1: D_{n-1} >M_{\rm{Tong}}&\approx&4.04.\\
a_{n}=2 \textrm{ and } a_{n+1}=2: D_{n-1} >M_{\rm{Tong}}&\approx&7.48.\\
a_{n}=2 \textrm{ and } a_{n+1}=3: D_{n-1} >M_{\rm{Tong}}&\approx& 10.92.
\end{eqnarray*}
If $a_{n}=1$ and $a_{n+1}=3$, then $F\approx 0.70$ and $G\approx 0.25$ So $\displaystyle{r-a_{n} >G}$ and $\displaystyle{\frac{1}{a_{n+1}}< R-a_{n+1} <F}$. Thus  
$$
D_{n-1}> \frac{a_{n+1}+1}{R-a_{n+1}} \approx  6.67 > M_{\rm{Tong}} \approx 5.76.
$$
For all other values of $a_{n}$ and $a_{n+1}$ either $D_{n-2}>r$ or $D_{n}>R$, or both.

\end{Ex}

\section{The case $D_{n-2}>r$ and $D_n>R$}
\label{sec:larger}
In this section we study the case that $D_{n-2}$ and $D_n$ are larger than given reals $r$ and $R$, respectively.
\begin{Theorem}
Let $r,R>1$ be reals, let $n\geq 1$ be an integer and let $F$ and $G$ be as given in~{\rm(\ref{def: F and G})}. 

If $D_{n-2}>r$ and $D_n>R$, then 
\begin{enumerate}
  \item if $\displaystyle{r-a_{n} \geq G}$ and $\displaystyle{R-a_{n+1} <F}$, then  
\[
D_{n-1}< \frac{a_{n+1}+1}{F},
\]
 \item if $\displaystyle{ r-a_{n} < G \textrm{ and } R-a_{n+1} \geq F}$, 
then  
\[
D_{n-1}< \frac{a_n+1}{G},
\]
\  \item if $r-a_{n}<\frac{1}{a_{n+1}+1}$ and $R-a_{n+1}<\frac{1}{a_{n}+1}$, then
\[
D_{n-1} < (a_{n}+1)(a_{n+1}+1),
\]
 \item in all other cases
\[
D_{n-1}<M_{\rm{Tong}}.
\]
\end{enumerate}
The bounds are sharp.
\label{th: larger}
\end{Theorem}
\begin{proof}
The proof is very similar to that of Theorem~\ref{th: intro smaller} The only `new' case is the one where  $r-a<\frac{1}{b+1}$ and $R-b<\frac{1}{a+1}$; see Figure~\ref{im: positions} ($v$).  If  $r-a<\frac{1}{b+1}$, then the graph of $f_{a,r}$ lies below $\Delta_{a,b} \subset \Omega$. Similarly, if $R-b<\frac{1}{a+1}$ the graph  $g_{b,R}$ lies left left of   $\Delta_{a,b} \subset \Omega$. In this case we have that $D_{n-2}>r$ and $D_{n}>R$ for all $(t_{n},v_{n}) \in \Delta_{a,b}$. In this case $D_{n-1}$ attains its maximum in the lower left corner $\left(\frac{1}{b+1},\frac{1}{a+1}\right)$. For the intersection point $(S,f_{a,r}(S))$ either $S<\frac{1}{b+1}$ or $f_{a,r}(S)<\frac{1}{a+1}$ and from Lemma~\ref{lem: increasing decreasing} we conclude $(a+1)(b+1)<M_{\rm Tong}$.
\end{proof}


\begin{Ex}
We again use  $r=2.9$ and $R=3.6$; see Figure~\ref{im: example largerr} and Table~\ref{tab: example larger}.

\begin{figure}[!ht]
\includegraphics[height=70mm]{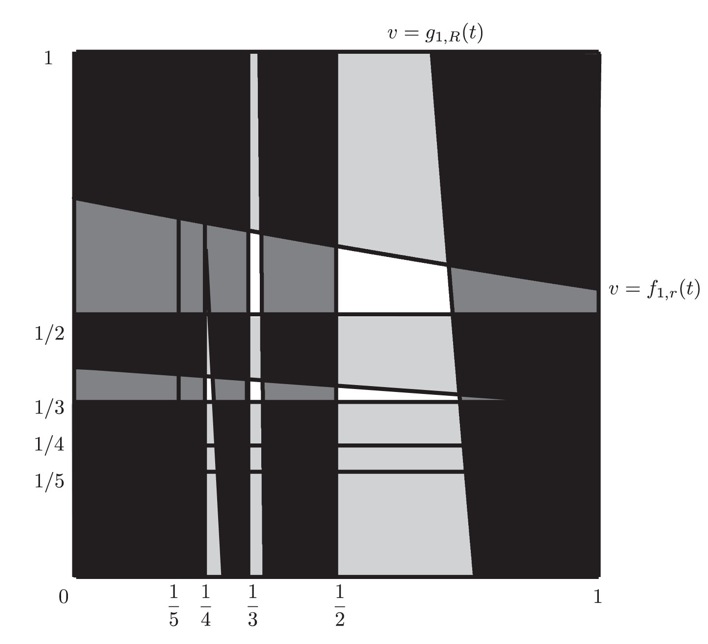}
\caption{Example with $r=2.9$ and $R=3.6$. The regions where $D_{n-2}>2.9$ are light grey, the regions where $D_n>3.6$ are dark grey. The intersection where both $D_{n-2}>2.9$ and $D_n>3.6$ is black.}
\label{im: example largerr}
\end{figure}

\begin{table}[!ht]
\begin{displaymath}
\begin{array}{ll|l|rr}
a_{n}&a_{n+1}&\textrm{Case}&\textrm{Upper bound for }D_{n-1} &\textrm{Tong's upper bound}\\
\hline
1 & 1 & (vi_a) & 2.30 & 2.30\\
1 & 2 & (vi_a) & 4.04 & 4.04\\
1 & 3 & (i) & 5.72 & 5.76\\
1 & 4 & (i) &  7.07 & 7.48\\
1 & 5,6,\dots & (i) &\dots &\dots\\
1 & 37 & (i) & 51.44  & 64.20\\
2 & 1 & (vi_c) & 4.04 & 4.04\\
2 & 2 & (vi_a) & 7.48 & 7.48\\
2 & 3 & (vi_a) & 10.92 & 10.92\\
2 & 4 & (i) &  13.79 & 14.36\\
2 & 5,6,\dots & (i) &\dots&\dots\\
2 & 42 & (i) & 116.00 & 144.97\\  
3 & 1 & (iii) & 4.04 & 5.76\\
3 & 2 & (iii) & 7.48 & 10.92\\
3 & 3 & (iii) & 10.92 & 16.08\\
3 & 4 & (v) &  13.79 & 21.23\\
4,5,6\dots & 1,2,3 & (iii) &\dots &\dots\\
3,4,5\dots & 4,5,6,\dots & (v)&\dots&\dots\\ 
17 & 29  & (v) & 540.00& 847.79 
\end{array}
\end{displaymath}
\caption{The sharp upper bounds and the Tong bounds for $D_{n-1}$ for $r=2.9$ and $R=3.6$. See Figure~\ref{im: positions} for cases ($i$)-($v$) and Figure~\ref{fig: msix} for ($vi_a$) and ($vi_c$).
}
\label{tab: example larger}
\end{table}

\end{Ex}

\newpage
\section{Asymptotic frequencies}
\label{sec: probabilities}

Due to Theorem~\ref{th: ergodic} and the ergodic theorem, the asymptotic frequency that an event occurs is equal to the measure of the area of this event in the natural extension. We calculate the measure of the region where $D_{n-2}>r$ and $D_n>R$. The same calculations can be done in the easier case where  $D_{n-2}<r$ and $D_n<R$.

\subsection{The measure of the region where $D_{n-2}>r$ and $D_n>R$ in a rectangle \texorpdfstring{$\Delta_{a,b}$}{} }
We calculate the measure in $\Delta_{a,b}$ above the graphs of $f_{a,r}$ and $g_{b,R}$ in the six cases from Figure~\ref{im: positions}. We denote $\log 2$ times the measure for case $(*)$ in  $\Delta_{a,b}$  by $m^{(*)}_{a,b}$.
\begin{eqnarray*}
m^{(i)}_{a,b} &=& \int_{\frac{1}{b+1}}^{\frac{1}{b}} \int_{f_{a,r}(t)}^{\frac{1}{a}} \frac{{\rm d}v \, {\rm d}t }{(1+tv)^{2}}  
= \int_{\frac{1}{b+1}}^{\frac{1}{b}}  \left[ \frac{-1}{t} \frac{1}{1+tv} \right]_{\frac{r}{a(r+1)+t}}^{\frac{1}{a}}     \,{\rm d}t\\
&=& \int_{\frac{1}{b+1}}^{\frac{1}{b}}  \frac{-1}{t} \frac{a}{a+t} + \frac{1}{t} \frac{a(r+1)+t}{(a+t)(r+1)} \,{\rm d}t\\
&=&  \int_{\frac{1}{b+1}}^{\frac{1}{b}}  \frac{-1}{t} +  \frac{1}{a+t} + \frac{1}{t}- \frac{r}{(a+t)(r+1)} \,{\rm d}t \\
&=&    \int_{\frac{1}{b+1}}^{\frac{1}{b}}   \frac{1}{(a+t)(r+1)} \,dt =  \frac{1}{(r+1)} \big[ \log(a+t) \big]_{\frac{1}{b+1}}^{\frac{1}{b}}\\
&=& \frac{1}{(r+1)}  \log\frac{(ab+1)(b+1)}{(ab+a+1)b}.
\end{eqnarray*}

Next we compute $m_{a,b}^{(v)}$, because it is handy for finding $m_{a,b}^{(ii)}$.
$$
m^{(v)}_{a,b}  = \int_{\frac{1}{b+1}}^{\frac{1}{b}} \int_{\frac{1}{a+1}}^{\frac{1}{a}} \frac{{\rm d}v \, {\rm d}t }{(1+tv)^{2}}  = \log \frac{(ab+1)(ab+a+b+2)}{(ab+a+1)(ab+b+1)}.\\
$$
For $m^{(ii)}_{a,b}$ we subtract the measure of the region in $\Delta_{a,b}$ below the graph of $f_{a,r}$ from $m_{a,b}^{(v)}$. 
$$
\begin{aligned}
m^{(ii)}_{a,b} &= m^{(v)}_{a,b}  - \int_{\frac{1}{b+1}}^{r-a} \int_{\frac{1}{a+1}}^{f_{a,r}(t)} \frac{{\rm d}v \, {\rm d}t }{(1+tv)^{2}}\\
&=   \log \frac{(ab+1)(ab+a+b+2)}{(ab+b+1)(ab+a+1)} - \frac{r}{r+1} \log \frac{r(b+1)}{ab+a+1}-\log \frac{ab+a+b+2}{(b+1)(r+1)}\\
&= \log \frac{(ab+1)(b+1)(r+1)}{(ab+b+1)(ab+a+1)} - \frac{r}{r+1} \log \frac{r(b+1)}{ab+a+1}.
\end{aligned}
$$
In the computation of $m^{(iii)}_{a,b}$ we use that $v=g_{b,R}(t)$  if and only if  $t = \frac{R}{v+b(R+1)}$, so
$$
m^{(iii)}_{a,b} =  \int_{\frac{1}{a+1}}^{\frac{1}{a}}   \int_{\frac{R}{b(R+1)+v}}^{\frac{1}{b}}  \frac{ {\rm d}t \, {\rm d}v }{(1+tv)^{2}}  
= \frac{1}{(R+1)}  \log\frac{(ab+1)(a+1)}{(ab+b+1)a}.
$$
Note that $m^{(iii)}_{a,b}$ is $m^{(i)}_{a,b}$ with $a$ interchanged with $b$ and $r$ replaced by $R$.

For $m_{a,b}^{(iv)}$ we find using the same techniques as before
$$
\begin{aligned}
m_{a,b}^{(iv)} &= m^{(v)}_{a,b}  -   \int_{\frac{1}{a+1}}^{R-b} \int_{\frac{1}{b+1}}^{\frac{R}{b(R+1)+v}} \frac{ {\rm d}t\, {\rm d}v  }{(1+tv)^{2}}\\
&=  \log \frac{(ab+1)(a+1)(R+1)}{(ab+a+1)(ab+b+1)} - \frac{R}{R+1} \log \frac{R(a+1)}{ab+b+1},
\end{aligned}
$$
which is  $m^{(ii)}_{a,b}$ where $a$ is interchanged with $b$ and $r$ replaced by $R$.

\begin{figure}[!ht]
\includegraphics[width=40mm]{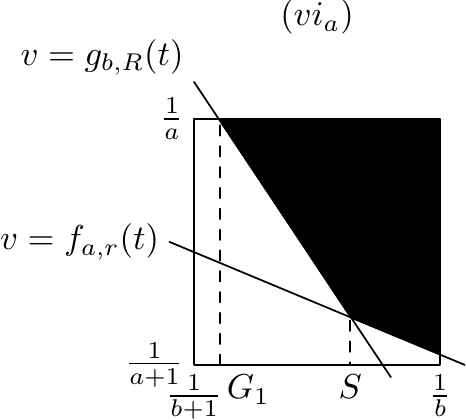}
\includegraphics[width=40mm]{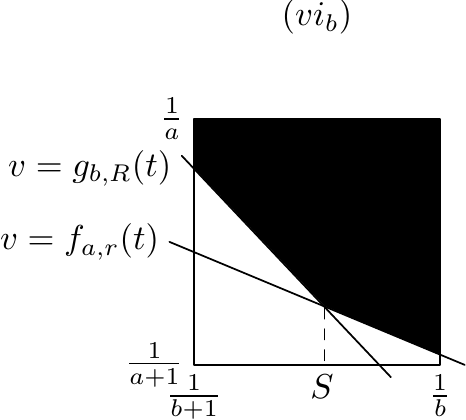}\\
\vskip5mm
\includegraphics[width=40mm]{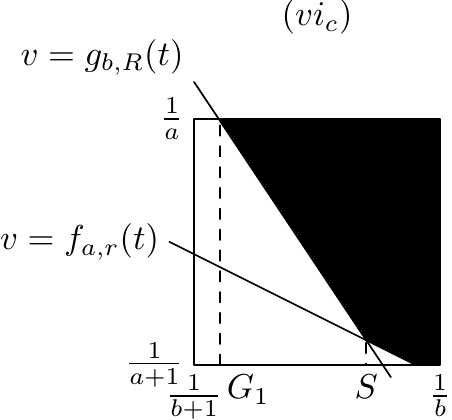}
\includegraphics[width=40mm]{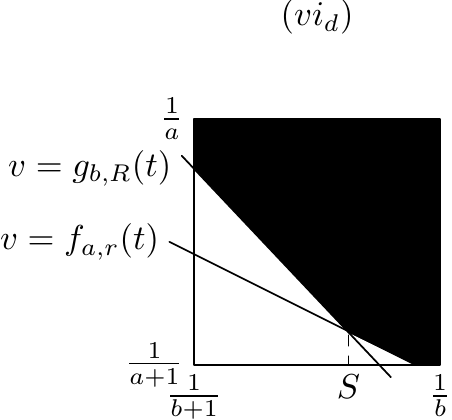}
\caption{The four possible configurations for case ($vi$).}
\label{fig: msix}
\end{figure}

In case ($vi$) there are four possibilities for the measure of the part above the graphs of $f_{a,r}$ and $g_{b,R}$, depending on where the graphs intersect with $\Delta_{a,b}$; see Figure~\ref{fig: msix}.

Denote $G_{1} = \frac{Ra}{ab(R+1)+1}$ (found from solving $g_{b,R}(G_{1})=\frac{1}{a})$ and recall from Lemma~\ref{lem: intersection f and g} that $S$ is the first coordinate of the intersection point of the graphs of $f_{a,r}$ and $g_{b,R}$. In this case we have that $(S,f_{a,r}(S) \in \Delta_{a,b}$.

\begin{itemize}
\item[ {(}\emph{vi}$_a${)}] If $r-a \geq \frac{1}{b}$ and $R-b \geq \frac{1}{a}$, then
\[
m^{(vi_a)}_{a,b} = \int_{G_1}^{S} \int_{g_{b,R}(t)}^\frac{1}{a}  \frac{{\rm d}v \, {\rm d}t }{(1+tv)^{2}}  + \int_{S}^{\frac{1}{b}} \int_{f_{a,r}(t)}^\frac{1}{a}  \frac{{\rm d}v \, {\rm d}t }{(1+tv)^{2}}.
\]
\item[{(}\emph{vi}$_b${)}] If $r-a \geq \frac{1}{b}$ and $R-b < \frac{1}{a}$, then
\[
m^{(vi_b)}_{a,b} = \int_{\frac{1}{b+1}}^{S} \int_{g_{b,R}(t)}^\frac{1}{a}  \frac{{\rm d}v \, {\rm d}t }{(1+tv)^{2}}  + \int_{S}^{\frac{1}{b}} \int_{f_{a,r}(t)}^\frac{1}{a}  \frac{{\rm d}v \, {\rm d}t }{(1+tv)^{2}}.
\]
\item[{(}\emph{vi}$_c${)}] If $r-a < \frac{1}{b}$ and $R-b \geq \frac{1}{a}$, then
\[
m^{(vi_c)}_{a,b} = \int_{G_1}^{S} \int_{g_{b,R}(t)}^\frac{1}{a}  \frac{{\rm d}v \, {\rm d}t }{(1+tv)^{2}}  + \int_{S}^{r-a} \int_{f_{a,r}(t)}^\frac{1}{a}  \frac{{\rm d}v \, {\rm d}t }{(1+tv)^{2}} + \int_{r-a}^{\frac{1}{b}} \int_{\frac{1}{a+1}}^\frac{1}{a}  \frac{{\rm d}v \, {\rm d}t }{(1+tv)^{2}}.
\]
\item[ {(}\emph{vi}$_d${)}] If $r-a < \frac{1}{b}$ and $R-b < \frac{1}{a}$, then
\[
m^{(vi_d)}_{a,b} = \int_{\frac{1}{b+1}}^{S} \int_{g_{b,R}(t)}^\frac{1}{a}  \frac{{\rm d}v \, {\rm d}t }{(1+tv)^{2}}  + \int_{S}^{r-a} \int_{f_{a,r}(t)}^\frac{1}{a}  \frac{{\rm d}v \, {\rm d}t }{(1+tv)^{2}} + \int_{r-a}^{\frac{1}{b}} \int_{\frac{1}{a+1}}^\frac{1}{a}  \frac{{\rm d}v \, {\rm d}t }{(1+tv)^{2}}.
\]
\end{itemize}

Using the following intergrals
$$
\begin{aligned}
 \int_{x}^{S} \int_{g_{b,R}(t)}^\frac{1}{a}  \frac{{\rm d}v \, {\rm d}t }{(1+tv)^{2}}  &=   \frac{1}{R+1} \log \frac{S(1-bx)}{x(1-bS)}+\log \frac{x(S+a)}{S(x+a)},\\
\int_{S}^{y} \int_{f_{a,r}(t)}^\frac{1}{a}  \frac{{\rm d}v \, {\rm d}t }{(1+tv)^{2}} &=\frac{1}{r+1}\log \frac{a+y}{a+S},\\
\int_{r-a}^{\frac{1}{b}} \int_{\frac{1}{a+1}}^\frac{1}{a}  \frac{{\rm d}v \, {\rm d}t }{(1+tv)^{2}} &= \log \frac{(ab+1)(r+1)}{(ab+b+1)r},
\end{aligned}
$$
we find that
$$
\begin{aligned}
m^{(vi_a)}_{a,b} & = \frac{1}{R+1} \log \frac{S(1-bG_1)}{G_1(1-bS)}&+& \frac{1}{r+1}\log \frac{ab+1}{(a+S)b}& +&  \log \frac{G_1(S+a)}{S(G_1+a)},\\
m^{(vi_b)}_{a,b} &=    \frac{1}{R+1} \log \frac{S}{(1-bS)} &+& \frac{1}{r+1}\log \frac{ab+1}{(a+S)b}  &+& \log \frac{S+a}{S(ab+a+1)},\\
m^{(vi_c)}_{a,b} 
&=\frac{1}{R+1} \log \frac{S(1-bG_1)}{G_1(1-bS)}& +& \frac{1}{r+1}\log \frac{r}{a+S}& +& \log \frac{G_1(S+a)(ab+1)(r+1)}{S(G_1+a)(ab+b+1)r},\\
m^{(vi_d)}_{a,b} &= \frac{1}{R+1} \log \frac{S}{(1-bS)}  &+&  \frac{1}{r+1}\log \frac{r}{a+S}& +&   \log \frac{(S+a)(ab+1)(r+1)}{S(ab+a+1)(ab+b+1)r} .
\end{aligned}
$$

\subsection{The total measure of the region where $D_{n-2}>r$ and $D_n>R$ in \texorpdfstring{$\Omega$}{the natural extension} }
For every $r>1$ and $R>1$ the asymptotic frequency that $D_{n-2}>r$ and $D_{n}>R$ can be found by adding a finite number of integrals. Let $\{x\} = x - \lfloor x \rfloor$ and $1_{A}$ be the indicator function of $A$, i.e.
\[
1_{A} = \left\{ \begin{array}{ll}
1 & \textrm{ if condition } A \textrm{ is satisfied,}\\ 
0 & \textrm{ else.}
 \end{array}\right.
\]

\begin{Proposition}
For almost all $x \in [0,1)$, and for all $r,R \geq 1$, we have that
\[
\log 2 \lim_{n \rightarrow \infty} \frac{1}{n} \#  \left\{ 2 \leq j \leq n+1;D_{j-2} > r \textrm{ and } D_{j}>R \right\} 
\]
exists and equals
\begin{eqnarray}
\sum_{a=1}^{\lfloor r \rfloor -1} \sum_{b=\lfloor R \rfloor +1}^{\infty} m^{(i)}_{a,b}   &+&  \sum_{a=1}^{\lfloor r \rfloor -1} \left(1_{(\{R\} \leq F)} m^{(i)}_{a,\lfloor R \rfloor} +  1_{(\{R\} \geq \frac1a)} m^{(vi_a)}_{a,\lfloor R \rfloor} + 1_{(F < \{R\} < \frac1a)} m^{(vi_b)}_{a,\lfloor R \rfloor} \right) \nonumber  \\
 + \sum_{a=1}^{\lfloor r \rfloor -1} \sum_{b=1}^{\lfloor R \rfloor -1} m^{(vi_a)}_{a,b}  
 &+& \sum_{b=\lfloor R \rfloor +1}^{\infty} \left( 1_{(\{r\} \geq \frac1b)} m^{(i)}_{\lfloor r \rfloor,b} +1_{(\frac{1}{b+1} <\{r\}<\frac1b)} m^{(ii)}_{\lfloor r \rfloor,b} \right) \nonumber \\
+  M_{r,R} &+& \sum_{b=1}^{\lfloor R \rfloor -1} \left(1_{(\{r\} \leq G)} m^{(iii)}_{\lfloor r \rfloor,b} +  1_{(\{r\} \geq \frac1b)} m^{(vi_a)}_{\lfloor r \rfloor,b} +1_{(G <\{r\}<\frac1b)} m^{(vi_c)}_{\lfloor r \rfloor,b} \right) \nonumber \\
 +  \sum_{a=\lfloor r \rfloor +1}^{\infty} \sum_{b=\lfloor R \rfloor +1}^{\infty} m^{(v)}_{a,b} &+&  \sum_{a=\lfloor r \rfloor +1}^{\infty} 
\left(1_{(\{R\} \geq \frac1a)} m^{(iii)}_{a,\lfloor R \rfloor} + 1_{(\{R\} \geq \frac1a)} m^{(vi_a)}_{a,\lfloor R \rfloor}+   1_{(F < \{R\}< \frac1a)} m^{(vi_b)}_{a,\lfloor R \rfloor}  \right) \nonumber \\
 + \sum_{a=\lfloor r \rfloor +1}^{\infty} \sum_{b=1}^{\lfloor R \rfloor -1}m^{(iii)}_{a,b}, && \nonumber
\end{eqnarray}
where $M_{r,R}$ is the measure of the regions where $D_{n-2}>r$ and $D_{n}>R$ in $\Delta_{\lfloor r \rfloor,\lfloor R \rfloor}$.

\end{Proposition}
\begin{proof}
Let $a,b \geq 1$ be integers. We denote strips with constant $a_n$ or $a_{n+1}$ by
$$
H_a=[0,1] \times \left[\frac{1}{a+1},\frac{1}{a}\right]  \quad \textrm{and} \quad V_b=\left [\frac{1}{b+1},\frac{1}{b}\right]\times [0,1].
$$

For $a < \lfloor r\rfloor$ the curve $v=f_{a,r}(t)$ is entirely inside the rectangle $H_a$ and (depending on the position of the curve $v=g_{b,R}(t)$) we are either in case ($i$) or ($vi$); see Figure~\ref{im: positions} and Remark~\ref{rmk: strips}. If $a > \lfloor r\rfloor $  the curve  $v=f_{a,r}(t)$ is entirely underneath $H_a$ and we are in case ($iii$), ($iv$) or ($v$).  For $a=\lfloor r\rfloor$ the curve $v=f_{a,r}(t)$ is partially inside and partially underneath $H_{\lfloor r \rfloor}$. In this strip we can have  each of the six cases.

Similarly, for $b< \lfloor R\rfloor$, the curve of $v=g_{b,R}(t)$ is entirely inside the rectangle $V_b$ and (depending on the position of the curve $v=g_{b,R}(t)$)  we are in case ($iii$) or ($vi$). For $b > \lfloor R \rfloor$  the curve  $v=g_{b,R}(t)$ is left of $V_b$ and we are in case ($i$), ($ii$) or ($v$) . For $b=\lfloor R\rfloor$ the curve $v=g_{b,R}(t)$ is partially inside and partially left of $V_{\lfloor R \rfloor}$ and we can have each of the six cases.

We use the strips $H_{ \lfloor r \rfloor} $ and  $ V_{\lfloor R \rfloor}$ to divide $\Omega$ in nine rectangles. Each of the nine terms in the sum in the proposition gives the measure of the region where $D_{n-2}>r$ and $D_{n}>R$  on one of those rectangles, we work from left to right and from top to bottom. The results follow from (\ref{eq: limit prob}), Remark~\ref{rmk: strips}, Theorem~\ref{th: larger} and the above. For instance, the first rectangle is given by $\left[0,\frac{1}{\lfloor R+1 \rfloor} \right) \times   \left[\frac{1}{\lfloor r \rfloor} ,1 \right)$ and we see that for every $\Delta_{a,b}$ in this rectangle we are in case ($i$). 
\end{proof}

\begin{Rmk}
All the infinite sums are just finite integrals, for example
\begin{equation}
\label{eq: left part}
\sum_{a=1}^{\lfloor r \rfloor -1} \sum_{b=\lfloor R \rfloor +1}^{\infty} m^{(i)}_{a,b}   = \int_{0}^{\frac{1}{\lfloor R \rfloor +1}} \int_{f_{a,r}(t)}^{\frac{1}{a}} \frac{{\rm d}v \, {\rm d}t }{(1+tv)^{2}}  .
\end{equation}
\end{Rmk}
%
\begin{Ex}
In this example we compute the asymptotic frequency that simultaneously $D_{n-2}>2.9$ and $D_{n}>3.6$; see Figure~\ref{im: example largerr} and Table~\ref{tab: probabilities}. Also compare with~Table~\ref{tab: example larger} where some of the upper bounds for this case are listed.
\begin{table}[!ht]
\begin{displaymath}
\begin{array}{rr|l|l}
a_{n}&a_{n+1}&\textrm{Case}&\textrm{asymptotic frequency}\\
\hline
1 & 1 & (vi_a) & 0.047\\
1 & 2 & (vi_a) & 0.025\\
1 & >2 & (i) & 0.106\\
2 & 1 & (vi_c) & 0.025 \\
2 & 2 & (vi_a) &0.013\\
2 & 3 & (vi_a) &0.090 \\
2 & >3 &  (i) &0.044  \\
>2 & 1 & (iii)& 0.097\\
>2 & 2 & (iii) & 0.050 \\
>2 & 3 & (iii)& 0.034 \\
>2 & > 3 & (v) & 0.115
\end{array}
\end{displaymath}
\caption{The probabilities that $D_{n-2}>2.9$ and $D_{n}>3.6$ in the various cases. 
}
\label{tab: probabilities}
\end{table} 

Summing over the cases yields that for almost all $x \in [0,1)\setminus \Q$ the asymptotic frequency that simultaneously $D_{n-2}>2.9$ and $D_{n}>3.6$ is 0.64.

We can also compute the conditional probability that $M_{\rm{Tong}}$ is the sharp bound. Given that  $D_{n-2}>2.9$ and $D_n>3.6$ the conditional probability that $M_{\rm{Tong}}$ is the sharp bound is 0.31.
\end{Ex}

\section{Results for $C_n$.}
\label{sec: C}
In~\cite{Tca04}, Tong states the following result as theorem without a proof.

\emph{
Let $t>1$, $T>1$ be two real numbers and
\[ K=\frac12 \left(  \frac{1}{t-1}+\frac{1}{T-1}+a_{n}a_{n+1}t\,T+\sqrt{\left(\frac{1}{t-1}+\frac{1}{T-1}+a_{n}a_{n+1}tT\right)^{2}-\frac{4}{(t-1)(T-1)}} \right).\]
Then
\begin{enumerate}
\item $C_{n-2}<t, C_{n}<T$ imply $C_{n-1}>K$;
\item $C_{n-2}>t, C_{n}>T$ imply $C_{n-1}<K$.
\end{enumerate}
}
This statement is not correct; assume for instance that $C_{n-2}<1.1$ and $C_{n}<1.4$, and that $a_{n}=a_{n+1}=1$. Part (1) of Tong's result then implies that \mbox{$C_{n-1}>11.94$}. However, by definition $C_{n-1} \in (1,2)$, so this bound is clearly wrong. 

In this section we give the correct result. The bounds in our theorems are sharp. We start with the case that both $C_{n-2}$ and $C_{n}$ are larger than given reals, this is related to the case where $D_{n-2}$ and $D_{n}$ are smaller than given numbers.

\begin{Theorem}
\label{th: c>t}
Let $t,T \in (1,2)$ and put
\begin{eqnarray*}
&&F' = \frac{a_{n+1}+1}{(a_{n}a_{n+1}+a_{n}+1)t-1},  \quad G' = \frac{a_{n}+1}{(a_{n}a_{n+1}+a_{n+1}+1)T-1} \\
\textrm{ and } && L'=t+T+a_{n}a_{n+1}tT-2.
\end{eqnarray*}
If $C_{n-2}>t$ and $C_{n}>T$, then
\begin{enumerate}
\item if $\displaystyle{\frac{1}{t-1}-a_{n} \geq G'}$ and  $\displaystyle{\frac{1}{T-1} - a_{n+1} < F' }$, then
\[
C_{n-1}< \frac{T}{(a_{n+1}+1)(T-1)}, 
\] 
\item if $\displaystyle{ \frac{1}{t-1}-a_{n}<G'}$ and  $\displaystyle{\frac{1}{T-1} - a_{n+1} \geq F' }$, then
\[
C_{n-1}< \frac{t}{(a_{n}+1)(t-1)}, 
\] 
\item in all other cases
\[
C_{n-1} < 1+\frac{L'-\sqrt{L'^{2}-4(t-1)(T-1)}}{2(t-1)(T-1)}.
\]
\end{enumerate}
The bounds are sharp.
\end{Theorem}
\begin{proof}
The proof follows from the fact that $C_{n}=1+\frac{1}{D_{n}}$ and Theorem~\ref{th: intro smaller}.  If $C_{n-2}>t$, then $D_{n-2} = \frac{1}{C_{n-2}-1}< \frac{1}{t-1}$ and likewise if $C_{n}>T$, then $D_{n}<\frac{1}{T-1}$. Setting $r=\frac{1}{t-1}$ and $R=\frac{1}{T-1}$, it directly follows from (\ref{def: F and G}) that $F=F'$ and $G=G'$. 

Consider case (1). The condition $\frac{1}{t-1}-a_n \geq G'$ is equivalent to $r-a_n\geq G$ and  $\displaystyle{\frac{1}{a_{n}+1} \leq  \frac{1}{T-1} - a_{n+1} < F' }$ is equivalent to $\displaystyle{\frac{1}{a_{n}+1} \leq R-a_{n+1} <F}$ in part (1) of Theorem~\ref{th: intro smaller}. We find that
$$
C_{n-1} < \frac{\frac{1}{T-1}-a_{n+1}}{a_{n+1}+1}+1 = \frac{T}{(a_{n+1}+1)(T-1)}.
$$
The proof of the second case is similar. For the third case we use Theorem~\ref{eq: MTong} for $M_{\rm{Tong}}$.
\begin{eqnarray*}
C_{n-1} &<& 1 + \frac{1}{M_{\rm{Tong}}} \\
&=& 1 + \frac{2}{t+T+a_{n}a_{n+1}tT-2+\sqrt{[t+T+a_{n}a_{n+1}tT-2]^{2}-4(t-1)(T-1)}}\\
&=& 1 + \frac{2}{L'+\sqrt{L'^{2}-4(t-1)(T-1)}} \cdot \frac{L'-\sqrt{L'^{2}-4(t-1)(T-1)}}{L'-\sqrt{L'^{2}-4(t-1)(T-1)}}\\
&=& 1 + \frac{L'-\sqrt{L'^{2}-4(t-1)(T-1)}}{2(t-1)(T-1)}.
\end{eqnarray*}
\end{proof}

\begin{Ex}
Take $t=1.1, T=1.4 $ and $a_{n}=a_{n+1}=1$. We find that $F'=0.870, G'=0.625$ and $L'=2.04$. Since $ \frac{1}{T-1}-a_{n+1}=\frac32 >F'$ case (1) of Theorem~\ref{th: c>t} does not apply. The second case does not apply either, since  $\frac{1}{t-1}-a_{n}=9 > G'$. So we are in case (3) and $C_{n-1}<1.50$.
\end{Ex}

We state the next theorem without a proof, since it is similar to that of Theorem~\ref{th: c>t}.  The only difference is that the proof is based on Theorem~\ref{th: larger} instead of Theorem~\ref{th: intro smaller}.
\begin{Theorem}
Let $t,T \in (1,2)$ and $F', G'$ and $L'$ be as defined in Theorem~\ref{th: c>t}. If $C_{n-2}<t$ and $C_{n}<T$, then
\begin{enumerate}
\item If $\displaystyle{\frac{1}{t-1}-a_{n} \geq G'}$ and  $\displaystyle{  \frac{1}{T-1} - a_{n+1} < F' }$, then
\[
C_{n-1}> 1+ \frac{F'}{a_{n+1}+1},
\] 
\item If $\displaystyle{G' \leq \frac{1}{t-1}-a_{n} }$ and $\displaystyle{ \frac{1}{T-1}-a_{n+1}<F' }$,  then
 \[
C_{n-1}>  1+ \frac{G'}{a_{n}+1}, 
\] 
\item If $\displaystyle{ \frac{1}{t-1}-a_{n}< \frac{1}{a_{n+1}+1} }$ and  $\displaystyle{ \frac{1}{T-1}-a_{n+1}< \frac{1}{a_{n}+1} }$, then
\[
C_{n-1} > 1 + \frac{1}{(a_{n}+1)(a_{n+1}+1)}.
\]
\item In all other cases
\[
C_{n-1} > 1+\frac{L'-\sqrt{L'^{2}-4(t-1)(T-1)}}{2(t-1)(T-1)}.
\]
\end{enumerate}
The bounds are sharp.
\end{Theorem}

\section*{Acknowledgements}
We thank Rob Tijdeman for his thoughtful comments on a previous version of this paper.

\bibliography{bibionica}
\bibliographystyle{abbrv}

\end{document}